
\def\addots{\mathinner{\mkern1mu\raise1pt\vbox{\kern7pt\hbox{.}}\mkern2mu
\raise4pt\hbox{.}\mkern2mu\raise7pt\hbox{.}\mkern1mu}}

\def\ad{\mathop{\text{\rm ad}}\nolimits}
\def\Ad{\mathop{\text{\rm Ad}}\nolimits}

\def\diag{\mathop{\text{\rm diag}}\nolimits}

%
%

\def\OO{\mathop{\text{\rm O}}\nolimits}

\def\Re{\mathop{\text{\rm Re}}\nolimits}

\def\Sl{\mathop{\text{\rm Sl}}\nolimits}
\def\SO{\mathop{\text{\rm SO}}\nolimits}
\def\SU{\mathop{\text{\rm SU}}\nolimits}

\def\Spec{\mathop{\text{\rm Spec}}\nolimits}

\def\tr{\mathop{\text{\rm tr}}\nolimits}

\def\0{\text{\bf 0}}
\def\1{\text{\bf {1}}}

\def\a{{\frak a}}

\def\c{\mathop{\text{\bf c}}\nolimits}

\def\e{{\frak e}}
\def\f{{\frak f}}
\def\g{{\frak g}}
\def\gl{{\frak {gl}}}
\def\h{{\frak h}}

\def\k{{\frak k}}
\def\l{{\frak l}}
\def\m{{\frak m}}

\def\n{{\frak n}}

\def\p{{\frak p}}
\def\q{{\frak q}}

\def\sp{{\frak {sp}}}

\def\su{{\frak {su}}}
\def\so{{\frak {so}}}
\def\sL{{\frak {sl}}}

\def\z{{\frak z}}

\def \res {\!\mid\!\!}

\def\C{{\Bbb C}} 
\def\D{{\Bbb D}} 
 
\def\N{{\Bbb N}}
\def\R{{\Bbb R}}

\def\:{\colon}  
\def\.{{\cdot}}
\def\|{\Vert}

\def \la {\langle}

\def \ra {\rangle}
\def\ssk{\smallskip}

\def\nin{\noindent}
\def\pagebreak{\vskip 0pt plus 0.0001fil\break}
\def\linebreak{\break}

\def\hat{\widehat}
\def\eps{\varepsilon}
\def\eset{\emptyset}
\def\nin{\noindent}
\def\oline{\overline}

\def\phi{\varphi}

\def\subeq{\subseteq}

\def\into{\hookrightarrow}

\input amstex
\documentstyle{gen-p}
\NoBlackBoxes

\issueinfo{00}
  {0}
  {Xxxx}
  {1997}

\topmatter
\title The c-function for non-compactly causal symmetric spaces and its  
relations to harmonic analysis and representation theory\endtitle

\rightheadtext{The c-function for NCC symmetric spaces}
\leftheadtext{Kr\"otz and \'Olafsson}

\author Bernhard Kr\"otz and Gestur \'Olafsson \endauthor

\address The Ohio State University, Department of Mathematics, 231 West 18th Avenue, 
Columbus, OH 43210-1174 \endaddress
\email kroetz\@math.ohio-state.edu \endemail

\thanks Supported in part by NSF-grant DMS-0097314.\endthanks

\address Louisiana State University, Department of Mathematics, Baton Rouge, LA 70803\endaddress
\email olafsson\@math.lsu.edu\endemail

\thanks Supported in part by NSF-grant DMS-0070607 and DMS-0139783.\endthanks

\subjclass 22E46 \endsubjclass

\endtopmatter
\document

\specialhead Introduction \endspecialhead 

The Harish-Chandra $\c$-function and its profound role in harmonic
analysis is well known to every one who has had any connections to
harmonic analysis on Riemannian symmetric spaces. The $\c$-function
shows up in a most natural way either as the coefficient of the leading
term in the asymptotic expansion of the spherical function; as the
density of the Plancherel measure; or
finally as the eigenvalue of the Knapp-Stein intertwining operator
acting on the one-dimensional space of $K$-fixed vectors in a principal series representation.
All of those different roles are related, but each one of them
explain different aspects of this jewel. The detailed knowledge
of the $\c$-function is given by the Gindikin-Karpelevic product
formula (cf.\ [5]) and the explicit representation of
the $\c$-function on rank one spaces in terms of $\Gamma$-functions.

\par The situation for non-Riemannian semisimple symmetric spaces
$G/H$ is quite different and much more involved. There are several
reason for this, the most obvious one being, that in general the
Hilbert space ${\Cal H}$ of an irreducible unitary representation $\pi$ of $G$ does in general
{\it not} have a nontrivial $H$-fixed vector. This facts leads us
into the study of $H$-invariant distribution vectors, i.e., understanding
the space $({\Cal H}^{-\infty})^H$. This space is finite dimensional, but in general the dimension
is greater than one. The action of an intertwining operator is therefore not
given by a scalar valued function but by a matrix. However, there is an important
class of symmetric spaces where the notion of spherical function and
{\it the} $\c$-function shows up naturally. Those are the
non-compactly causal symmetric spaces. A simple characterization
is in terms of an $H$-invariant open convex cone $C$ consisting of
hyperbolic elements and diffeomorphic to an open $H$-invariant subset
in $G/H$. This should be understood as a replacement for the Cartan decomposition
$G/K\simeq \p$ for Riemannian symmetric spaces. The cone $C$ 
is the domain where
the spherical functions $\phi_\lambda$ on $G/H$ are supported. The asymptotic behavior
of those functions is again given by a $\c$-function, which - again analogous to
the Riemannian case - can be represented as an integral over a unipotent
group $\overline N$. The new complication is, that the integral is now
not over the full group $\oline N$, but only over subset given by $\overline N\cap HP_{\text{\rm min}}$ where
$P_{\text{\rm min}}$ is certain minimal parabolic subgroup. The understanding
of this domain of integration is the first step in understanding the
$\c$-function and to generalize the Gindikin-Karpelevic product formula
to this setting.

\par In this overview article several aspects of the $\c$-function and its
role in harmonic analysis and representation theory are discussed. Generally we do not go into any details
of the proofs, but try instead to explain the ideas by our
main example $\Sl (n,\R )/\SO (p,q)$ and include references where detailed proofs
can be found. We also include references to articles which have provided important particular steps or
introduced new ideas and tools, although we have to point out that this list
is in no way complete and should not be taken as such.

\par In the first
three sections we introduce the minimal amount necessary for understanding
the fine structure of $\overline N\cap HP_{\text{\rm min}}$ that leads up to the product formula.
We start by introducing the category of non-compactly causal symmetric spaces, but
the most important part is the convex structure of the real bounded symmetric
space $\Omega \simeq H/H\cap K$.

\par The integral representation for the $\c$-function is introduced in Section 4.
As an example we give an explicit calculation for the simplest rank-one case \linebreak 
$\Sl (2,\R )/\SO_o(1,1)$. The result is that the $\c$-function
is given by a Beta-function. For the general rank-one non-compactly causal symmetric spaces
$\SO_o(n+1,1)/\SO_o(n,1)$ the formula is $\c (\lambda )=\kappa B(n/2, 1-(\lambda +n)/2)$
where $\kappa$ is a positive constant depending on the normalization of
measures.
We finish the section by the general product formula. Several important articles
leading up to the final proof in  [16] include the important special case
of Cayley type symmetric spaces worked out by J. Faraut in  [2] and
the work by P. Graczyk [6]
where the $\c$-function  for our main example
$\Sl (n,\R)/\SO (p,q)$ was calculated.

\par The second half of the article is devoted to the role of the
$\c$-function in harmonic analysis and representation theory.
We start by discussing the spherical functions $\phi_\lambda$
and their asymptotics.  We conclude Section 5 by the Harish-Chandra
type expansion formula for the spherical function where the
coefficients are now given by the $\c$-function for $G/H$. 
Finally we mention in Theorem 5.10 how the the numerator of
the $\c$-function gives information about the singularities of
$\phi_\lambda$ in the spectral parameter.

\par We discuss the $H$-spherical highest weight representations
of the $c$-dual group $G^c$ in Section
6. Those representations play an important role in harmonic analysis
on $G^c/H$ but here we focus on the space of $H$-spherical distribution
vectors. The main result is stated in Theorem 6.5, the Averaging Theorem.
Let $(\pi_\lambda, {\Cal H}_\lambda)$ be a unitary $H$-spherical
highest weight representation. We fix a normalized highest weight vector $v_\lambda\in {\Cal H}_\lambda$. 
Then for ``big'' parameters
we have that $\int_H\pi_\lambda (h)v_\lambda \, dh$ exists and equals to $\c (\lambda +\rho)\nu$
where $\nu$ is $H$-spherical and normalized by $\la\nu ,v_\lambda \ra=1$.
This results leads then to an explicit formula involving $\c$-function for the
{\it formal dimension}.

\par Section 5 and 6 are finally connected in Section 7 where we introduce the $H$-spherical
distribution $\theta_\lambda(f)= f\mapsto \la\pi_\lambda (f)\nu ,\nu\ra$. It turns out that this distribution
is in fact determined by the spherical function and the $\c$-function by representing
$\theta_\lambda$ as a hyperfunction. As such we have $\theta_\lambda = \c (\lambda + \rho)^{-1}\phi_{\lambda+\rho}$.

\head 1. Causal symmetric spaces \endhead 

In this section we recall some  basic facts and definitions that will be
needed later on. Our standard reference for the theory of {\it causal
symmetric spaces} is  is the monograph [10], in particular chapter 3.

\par Let $G$ be a connected semisimple Lie group with finite
center and $\g$ its Lie algebra. 
If $\sigma$ is an automorphism of $G$, then $\sigma$ will also
denote the derived automorphism $d\sigma(\1)$ of $\g$ (and the
complex linear extension to $\g_\C=\g\otimes_\R\C$) and vice versa. 
Let $\theta$ be a Cartan involution on $G$ and let $K=G^\theta=\left\{g\in G\, :\,
\theta(g)=g\right\}$ be  the corresponding maximal
compact subgroup. On the Lie algebra level $\theta$ induces the Cartan decomposition
$\g=\k + \p$ where $\k=\left\{X\in \g\, :\, \theta(X)=X\right\}$ is the Lie algebra of $K$
and $\p =\left\{X\in \g\, :\, \theta(X)=-X\right\}$ is the (-1)-eigenspace of $\theta$.
Let $\tau$ be a non-trivial involution on $\g$. We can and will assume that
$\tau$ commutes with $\theta$. Write $\g=\h + \q$ for the
$\tau$-eigenspace decomposition according to the eigenvalues $+1$ and $-1$.
As $\theta$ and $\tau$ commutes we get the joint eigenspace decomposition
$$\g=\k\cap \q +\k\cap \h + \p\cap \q + \p\cap \h \, .$$
Let  $H$ be an open subgroup of $G^\tau$ with $G^\tau$ the
$\tau$-fixed points of $G$. Then $G/H$ is called an {\it
affine symmetric space} or simply a {\it symmetric space}.
Notice that the Lie algebra of $H$ is $\h$.
\par In this article we will always assume that $G$ is contained in the
connected simply connected Lie group $G_\C$ with Lie algebra $\g_\C$. Furthermore
we will assume that $H=G^\tau$. Those two assumptions simplify some
arguments and our notation, but all the main results remain  valid
for the general case. Note in particular that every homomorphism $\sigma$ of $\g_\C$
corresponds to a homomorphism of $G_\C$ as $G_\C$ is simply connected.
\par The symmetric pair
$(\g ,\h)$ is called {\it irreducible} if the only
$\tau$-invariant ideals in $\g$ are the trivial ones,
$\left\{0\right\}$ and $\g$.  We say
that the symmetric space $G/H$ is {\it irreducible} if $(\g,\h)$ is
irreducible. If $(\g ,\h)$ is irreducible, then either $\g$ is
simple or $\g \simeq \h\oplus \h$, with $\h$ simple, and $\tau (X,Y)=
(Y,X)$.
In the second case we have that $G/H\ni (a,b)H\mapsto ab^{-1}\in H$
defines an isomorphism of $G/H$ onto $H$. From now on we will
always assume that $G/H$ is irreducible.

\definition{Definition 1.1} Let $C\not= \emptyset$ be an open convex subset of
$\q$. Then $C$ is called {\it hyperbolic}, if for all $X\in C$
the operator $\ad(X)$ is semisimple with real eigenvalues.
We say that $C$ is {\it elliptic}, if for all $X\in C$ the operator  $\ad(X)$
is semisimple  with imaginary spectrum.
\enddefinition

\definition {Definition 1.2} {\bf (Causal symmetric spaces}) Suppose that $G/H$ is irreducible. Then $G/H$ is
said to be {\it causal}
if there exists a non-empty  open $H$-invariant convex cone $C\subseteq \q$,
containing no affine lines.
\enddefinition

There are two different types of causal symmetric spaces,
the {\it non-compactly causal} symmetric spaces (NCC) and  the
{\it compactly
causal} symmetric spaces (CC). In addition, there is the intersection of
those two classes, the {\it Cayley type} symmetric spaces (CT).

\definition {Definition 1.3} {\bf (NCC)} Assume that $G/H$ is an irreducible causal symmetric space. Then
the following two condirtions are equivalent:
\roster 
\item  There exists a non-empty $H$-invariant open hyperbolic cone $C
\subseteq \q$ which contains  no affine lines;

\item There exists an element $T_0\in \q\cap \p$, $T_0\not= 0$,
which is fixed by $H\cap K$.
\endroster 
If one of those equivalent conditions are satisfied, then $G/H$ is called
{\it non-compactly causal}.
\enddefinition

\definition {Definition 1.4} {\bf (CC)} Assume that $G/H$ is an irreducible causal symmetric space. Then
the following two conditions are equivalent:
\roster 
\item  There exists a non-empty  $H$-invariant open elliptic cone $C
\subseteq \q$ which contains no affine lines;

\item There exists an element $X_0\in \q\cap \k$, $X_0\not= 0$,
which is fixed by $H\cap K$.
\endroster 
If one of those equivalent conditions are satisfied, then $G/H$ is
called {\it compactly causal}.
\enddefinition

\definition{Definition 1.5} {\bf (CT)} Assume that $G/H$ is an irreducible causal
symmetric space. Then $G/H$ is called a {\it symmetric space of
Cayley type}, if it is both non-compactly causal and compactly
causal.
\enddefinition

\remark{Remark 1.6} (a) The elements $T_0$ and $X_0$ in Definition 1.3
and Definition 1.4 are
unique up to multiplication by scalar. If $G/H$ is NCC then we can,
and will, normalize
$T_0$ such that $\ad (T_0)$ has spectrum $\left\{0,1,-1\right\}$. The eigenspace
corresponding to $0$ is exactly $\g^{\theta\tau}=\h\cap \k + \p\cap \q$.
If $G/H$ is
CC, then we normalize $X_0$ such that the spectrum of
$\ad(iX_0)$ is $\left\{0,1,-1\right\}$. In this case the zero eigenspace
is exactly $\k$.

\par\nin (b) If $G/H$ is compactly causal, then $\k$ has a non trivial
center $\z (\k )$ and $\z (\k )\cap \q=\R X_0$. If $G$ is simple then $\z (\k )
=\R X_0$. If $\g \simeq \h\oplus \h$. Then $X_0 = (X_0^\prime ,- X_0^\prime )$
with $X_0^\prime $ central in $\h\cap\k $ and
$\z (\k )=\R X_0\oplus \R (X_0^\prime ,X_0^\prime )$.
\endremark 

\par
Denote the complex linear extension of $\tau$ to $\g_\C$
by $\tau$.
The $c$-dual $(\g^c,\tau^c)$ of $(\g,\tau)$ is defined
by $\g^c=\h\oplus i\q\subeq \g_\C$ with involution
$\tau|_{\g^c}$. Notice that the $c$-dual of $(\g^c,\tau^c)$
is $(\g ,\tau)$. Let $G^c$ denote the analytic subgroup
of $G_\C$ with Lie algebra $\g^c$. As we are assuming that
$G_\C$ is simply connected and $H=G^\tau$ it follows that 
$G\cap G^c=H$. Hence $G^c/H$ is well defined. The
space $G^c/H$ is called the {\it c-dual} of $G/H$.

We recall the following fact:

\proclaim {Proposition 1.7} Assume that $G/H$ is an irreducible symmetric space.
Then the following holds:
\roster 
\item $G/H$ is non-compactly causal if and only if $G^c/H$ 
is compactly causal.

\item  $G/H$ is of Cayley type if and only if
$G/H\simeq G^c/H$.

\item Suppose that $G/H$ is causal. Then
$G/H$ is of Cayley type if and only if $H$ is not simple. In that
case $\z (\h )$ is one dimensional and contained in $\p$. 
\endroster 
\quad\qed 
\endproclaim 

\example{Example 1.8} {\bf (The group case)} Let $H$ be a connected semisimple Lie group,
$G=H\times H$, and $\tau (a,b)=(b,a)$. Then $G/H\simeq H$ is
compactly causal if and only if $H$ is a Hermitian group. In this case $G_\C=H_\C\times H_\C$ and
$G^c/H\simeq H_\C/H$ (cf.\ [10, Example 1.2.2]).
\endexample

We now come to our main example which will accompany us through the 
whole article. 

\example{Example 1.9} {\bf (Main example)} Let $G=\Sl(n,\R)$. Then $G_\C=\Sl(n,\C)$ is simply connected. 
We take the standard Cartan involution $g\mapsto \theta(g)=(g^{-1})^{t}$ for $g\in G$. 
Then $K=\SO(n,\R)$. 
Let $1<p<n$ and set $q=n-p$. Define an $n\times n$-matrix $I_{p,q}$ by 
$$I_{p,q}=\pmatrix \1_{p\times p}& 0\\
 0 &-\1_{q\times q}\endpmatrix
\ .$$

Then 
$$\tau\: G\to G, \ \ g\mapsto I_{p,q} \theta(g) I_{p,q}$$
defines an involution on $G$ commuting with $\theta$. We take $H=G^\tau=\SO(p,q)$. Notice that 
$$H\cap K=S(\OO(p)\times \OO(q))\ .$$
We have 
$$\q\cap\p=\left\{ \pmatrix  A& 0\\  0 &B\endpmatrix \: A\in M(p,\R),\  B\in M(q,\R);\ 
A^t=A,\  B^t=B,\  \tr (A)+\tr (B)=0\right\}\ .$$
Furthermore the element
$$T_0= {1\over n}\pmatrix q\1_{p\times p}& 0\\ 0 &-p\1_{q\times q}\endpmatrix \leqno(1.1)$$
is in $\q\cap\p$ and $T_0$ is $H\cap K$-fixed. According to Definition 1.3 the symmetric 
spaces $G/H=\Sl(n,\R)/ \SO(p,q)$ are NCC for all $1<p<n$. The dual space is readily 
computed to be 
$$G^c/H=\SU(p,q)/ \SO(p,q)\ .$$
We will continue the discussion of this example in the sequel. 
\endexample 

For the convenience to the reader we list here 
the three tables of causal symmetric pairs
(cf.\ [10, Th.\ 3.2.8]).
\bigskip
\centerline{\bf Table I}
\centerline{\bf The group case}
$$\vbox{\tabskip=0pt\offinterlineskip
\def\tablerule{\noalign{\hrule}}
\halign{\strut#&\vrule#\tabskip=1em plus2em&
\hfil#\hfil&\vrule\hfil#\hfil&\hfil#\hfil &\vrule\hfil#\hfil&
\hfil#\hfil&\vrule\hfil#\hfil \tabskip=0pt\cr\tablerule
&&\omit\hidewidth $\g$ \hidewidth&& 
\omit\hidewidth $\g^c$  \hidewidth&& 
\omit\hidewidth $\h$\hidewidth&  \cr
&&\omit\hidewidth CC\hidewidth
&& \omit\hidewidth NCC\hidewidth && &\cr\tablerule
&& $\su (p,q)\oplus \su(p,q)$ && $\sL (p+q,\C)$ && $\su (p,q)$
&\cr\tablerule
&& $\so^*(2n)\oplus \so^*(2n)$ && $\so (2n,\C)$ && $\so^*(2n)$ &
 \cr\tablerule
&& $\so (2,n)\oplus \so (2,n)$ && $\so(n+2,\C)$ && $\so (2,n)$ & 
\cr\tablerule
&& $\sp (n,\R )\oplus\sp (n,\R)$ && $\sp (n,\C)$ && $\sp (n,\R)$
&\cr\tablerule
&& $\e_{6(-14)}\oplus \e_{6(-14)}$ && $\e_{6}$ &&
$\e_{6(-14)}$ 
& \cr\tablerule
&& $\e_{7(-25)}\oplus \e_{7(-25)}$ && $\e_7$ && $\e_{7(-25)}$ &\cr\tablerule
}}
$$

\bigskip
\centerline{\bf Table II}
\centerline {\bf $\g$ simple and not CT}
$$\vbox{\tabskip=0pt\offinterlineskip
\def\tablerule{\noalign{\hrule}}
\halign{\strut#&\vrule#\tabskip=1em plus2em&
\hfil#\hfil&\vrule\hfil#\hfil&\hfil#\hfil &\vrule\hfil#\hfil&
\hfil#\hfil&\vrule\hfil#\hfil \tabskip=0pt\cr\tablerule
&&\omit\hidewidth $\g$ \hidewidth&& 
\omit\hidewidth $\g^c$  \hidewidth&& 
\omit\hidewidth $\h$\hidewidth&  \cr
&&\omit\hidewidth CC\hidewidth
&& \omit\hidewidth NCC\hidewidth && &\cr\tablerule
&& $\su (p,q)$ && $\sL (p+q,\R)$ && $\so (p,q)$
&\cr\tablerule
&& $\su^* (2p,2q)$ && $\su^*(2(p+q)$ && $\sp(p,q)$ &
 \cr\tablerule
&& $\so^*(2n)$ && $\so(n,n)$ && $\so (n,\C)$ & 
\cr\tablerule
&&$\so (2,n)$&&$\so (1,n+1)$&&$\so (1,n)$&\cr\tablerule
&& $\so(2,p+q)$ && $\so(p+1,q+1)$ && $\so(1,p)\oplus \so (1,q)\, ,\, p,q\ge 2$
&\cr\tablerule
&& $\sp (2n,\R)$ && $\sp (n,\R)$ &&
$\sp (n,\C)$ 
& \cr\tablerule
&& $\e_{6(-14)}$ && $\e_{6(6)}$ && $\sp(2,2)$ &\cr\tablerule
&&$\e_{6(-14)}$&&$\e_{6(-26)}$ &&$\f_{4(-20)}$&\cr\tablerule
&&$\e_{7(-25)}$&&$\e_{7(7)}$&&$\su^*(8)$&
\cr\tablerule
}}
$$

\bigskip

\centerline{\bf Table III} 
\centerline{\bf $\g$ simple and CT} 
$$\vbox{\tabskip=0pt\offinterlineskip
\def\tablerule{\noalign{\hrule}}
\halign{\strut#&\vrule#\tabskip=1em plus2em&
\hfil#\hfil&\vrule\hfil#\hfil&\hfil#\hfil &\vrule\hfil#\hfil&
\hfil#\hfil&\vrule\hfil#\hfil \tabskip=0pt\cr\tablerule
&&\omit\hidewidth $\g$ \hidewidth&& 
\omit\hidewidth $\g^c\simeq \g$  \hidewidth&& 
\omit\hidewidth $\h$\hidewidth&  \cr\tablerule
&&$\sp(n,\R)$&&$\sp (n,\R)$&&$\gl (n,\R)$&\cr\tablerule
&& $\su (n,n)$ && $\su (n,n)$ && $\sL (n,\C)\oplus \R$
&\cr\tablerule
&& $\so^*(4n)$ && $\so^*(4n)$ && $\su^*(2n)\oplus \R$ &
 \cr\tablerule
&&$\so(2,n)$&&$\so (2,n)$&&$\so(1,n-1)\oplus \R$&\cr\tablerule
&&$\e_{7(-25)}$&&$\e_{7(-25)}$&&$\e_{6(-26)}\oplus \R$&
\cr\tablerule
}}
$$

\head 2. Bounded symmetric domains \endhead 

In this section we discuss the relation between causal symmetric spaces
and totally real symmetric subdomains of bounded complex symmetric domains.
The $\c$-function, that we discuss later, turns out to be an integral
over those domains.
We begin by recalling some elementary structure theory for causal symmetric spaces
and introduce some notation that will be used later on. 

\par
{}From now on $G/H$ will denote a NCC symmetric space. Let $\a\subeq \p\cap \q$ 
be a maximal abelian subspace. Note that NCC implies that $\a$ is also maximal 
abelian in $\p$ and in $\q$. Let $\m=\z_\k(\a)$. As $\a$ is maximal abelian in $\q$ and $\p$  it follows that
$\m=\z_\h (\a)$. For $\alpha \in \a^* $ let $\g^\alpha =\left\{X\in \g\: (\forall H\in\a)\,
[H,X]=\alpha (H)X\right\}$ and let $\Delta =\left\{\alpha \in \a^*\: \g^\alpha\not= \left\{0\right\}\right\}\setminus \left\{0\right\}$
denote the set of restricted roots with respect to $\a$. 
Then we have the root space decomposition 
$$\g=\a\oplus\m\bigoplus_{\alpha\in\Delta}\g^\alpha$$
associated to $\Delta$.
\par Recall the $H\cap K$-fixed element $T_0\in \q\cap \p$, 
and notice that $T_0\in\a$ no matter what our choice of $\a\subeq \q\cap \p$ will be 
(cf.\ [10, Lemma 1.3.5]). As $\Spec(\ad T_0)=\left\{-1, 0, 1\right\}$
we obtain a triangular decomposition 
$$\Delta=\Delta_n^-\amalg \Delta_0\amalg \Delta_n^+\leqno(2.1)$$
with $\Delta_n^\pm=\left\{ \alpha\in \Delta\: \alpha(T_0)=\pm 1\right\}$ and 
$\Delta_0=\left\{\alpha\in \Delta\: \alpha(T_0)=0\right\}$.  
Using (2.1) we obtain a triangular decomposition 
of $\g$
$$\g=\n^-\oplus \g(0)\oplus \n^+\leqno(2.2)$$
with 
$$\n^-=\bigoplus_{\alpha\in \Delta_n^-}\g^\alpha,\quad \n^+=\bigoplus_{\alpha\in \Delta_n^+}\g^\alpha\quad
\hbox{and} \quad \g(0)=\a\oplus\m\bigoplus_{\alpha\in\Delta_0}\g^\alpha=\h\cap\k +\q\cap \p \ . $$
Notice that $\z(\g(0))=\R T_0$ (cf.\ [10, Lemma 1.3.5]). Denote by
$N^+$ respectively $N^-$ the analytic subgroup of $G$ corresponding to the
Lie algebra $\n^+$ respectively $\n^-$. Then both $N^+$ and $N^-$ are simply
connected abelian Lie groups. Let $G(0)=Z_G(T_0)$. Then the Lie algebra
of $G(0)$ is $\g(0)$. By our assumption that $G$ is
contained in a simply connected group $G_\C$ and
$H=G^\tau$ we get $G(0)=(H\cap K)\exp (\p\cap \q)=G\cap K^c_\C$. Notice
that by the same argument we also have
$M:=Z_K(\a )=Z_{H\cap K}(\a)\subseteq H\cap K$. Finally we set $P_{\text{\rm max}}=
G(0)N^+$. Then
$P_{\text{\rm max}}$ is a maximal parabolic subgroup of $G$.
\par Recall the Hermitian $c$-dual $\g^c=\h+ i\q$ of $\g$. Then $\k^c=
(\h\cap\k)+ i(\q\cap \p)$
is a maximal compact subalgebra of $\g^c$. Notice that $\k^c$ and $\g(0)$
have the same complexification. Let $\p^+:=(\n^+)_\C$ be the complexification of $\n^+$,
and similarily  $\p^-:=(\n^-)_\C$. Thus complexifying (2.2) yields
$$\g_\C=\p^-\oplus(\k^c)_\C \oplus\p^+, \leqno(2.3)$$
the familiar Harish-Chandra decomposition of $\g_\C$ associated
to its Hermitian real form $\g^c$.

\par Write $\sigma\:\g_\C\to \g_\C$ for the complex conjugation of $\g_\C$
with respect to the real form $\g$. Notice that $\sigma$ leaves the decomposition
(2.3) stable. Taking the space of $\sigma$-fixed points in (2.3) we get the decomposition in  (2.2) back.

\par Let $P^+$ and $P^-$ be the analytic subgroups of $G_\C$ corresponding 
to $\p^+$ respectively $\p^-$. Similarily we define $K^c$ and $K_\C^c$. By a
well known result of Harish-Chandra we have $G^c\subeq P^- K_\C^c P^+$. Hence
the image of the Hermitian space $G^c/K^c$ under the embedding 
$$\iota\: G^c/K^c\into G_\C/ K_\C^c P^+, \ \ gK^c\mapsto gK_\C^cP^+\leqno(2.4)$$
is contained in $P^-K_\C P^+$. As $P^-\cap K_\C^c P^+=\left\{\1\right\}$,  it
follows that we can realize $G^c/K^c$ in $P^-$ and hence, as
$\exp|_{\p^-} :\p^-\to P^-$ is a analytic diffeomorphism, also in $\p^-$.
We write ${\Cal D}=(\exp|_{\p^-})^{-1}(\iota (G^c/K^c))$.
Then ${\Cal D}\subeq \p^-$ is a bounded symmetric  domain. Notice that
$K^c_\C P^+\cap G= P_{\text{\rm max}}$. We get therefore a $G$-equivariant embedding
of flag manifolds
$$G/P_{\text{\rm max}}\into G^c_\C/K^c_\C P^+\ . $$
Furthermore
$$(G^cK_\C P^+)^{\sigma} =HP_{\text{\rm max}}\ .
\leqno(2.5)$$

\par By (2.4) it follows that $\sigma$ induces an anti-holomorphic involution on ${\Cal D}$.
Write $\Omega :=\left\{X\in {\Cal D}\: \sigma(X)=X\right\}\subseteq \n^-$.
Using that $\iota$ intertwines with $\sigma$, we obtain from (2.5) that
$$\Omega\simeq H/ H\cap K \leqno(2.6)$$
and
$$HP_{\text{\rm max}}=\exp (\Omega )P_{\text{\rm max}}\ . \leqno(2.7)$$
As ${\Cal D}$ is convex, the same then holds for $\Omega$. More precisely
we have (Herman's convexity theorem)

$$\Omega=\left\{ X\in \n^-\: \|\ad(X+\tau(X))\| <1\right\}\ .\leqno(2.8)$$
Here $\|\cdot\|$ stands
for the operator norm related to the inner product
$(X,Y)=-\tr (\ad (X)\ad (\theta(Y))$ on $\g$.

\example{Example 2.1} Assume that $G$ is complex. Then $G/H=H_\C/ H$ where $H$ is of
Hermitian type. In this case $\Omega={\Cal D}_H$ with ${\Cal D}_H$ the
Harish-Chandra realization of the Hermitian symmetric space $H/ H\cap K$.
\endexample

\example{Example 2.2} {\bf (Main example continued)} We continue our discussion
of the NCC space $G/H=\Sl(n,\R)/\SO(p,q)$ from Example 1.9. {}From now on we will
always assume that $p\leq q$.
\par Let $\a$ be the space of diagonal
matrices
$$\a=\left\{\diag(t_1,\ldots, t_n)\: t_j\in \R, \ \sum_{j=1}^n t_j=0\right\}\ .$$
Then $\a$ is maximal abelian in $\p\cap \q$ and it is clear
that $\a$ is also maximal
abelian in $\p$ and $\q$. The root system $\Delta$ is of
type $A_{n-1}$ and is given by
$\Delta =\left\{\eps_k-\eps_l\: k\neq l\right\}$ where
$\eps_j\in \a^*$ is defined by
$$\eps_j(\diag(t_1,\ldots, t_n))=t_j\ . $$
We choose $\Delta^+=\left\{\eps_i-\eps_j\: i<j\right\}$ as a positive
system. Let $E_{ij}=(\delta_{ki}\delta_{lj})\in M(n\times n,\R)$.
In this notation the root spaces are given by
$$\g^{\eps_{i}-\eps_{j}}=\R E_{ij}\, .$$
Furthermore by the definition of $T_0$ in  (1.1)  we obtain that
$$\Delta_n^+=\left\{\eps_k-\eps_l\: 1\leq k\leq p, \ p+1\leq l\leq n\right\}\ .$$

\par
Using that $\Delta_n^-=-\Delta_n^+$ and
$\n^-=\theta (\n^+)=\left\{X^t\:X\in \n^+\right\}$
we get:

$$\n^+=\left\{\pmatrix \0_{p\times p} & A\\ \0_{q\times p}& \0_{q\times q}\endpmatrix 
\: A\in M(p\times q,\R)\right\}\simeq M(q\times p,\R)
$$
and
$$\n^-=\left\{\pmatrix \0_{p\times p} & \0_{p\times q}\\
B & \0_{q\times q}\endpmatrix 
\: B\in M(q\times p,\R)\right\}\ .$$
Define elementary matrices in $M(q\times p, \R)$ by $F_{ij}=(\delta_{ki}\delta_{lj})$. 
Then we obtain a vector space isomorphism
$$\n^-\simeq M(q\times p, \R)\ .$$
by linear extension  of the basis vectors assignments $E_{i+p, j}\leftrightarrow F_{ij}$ for 
$1\leq i\leq q$ and $1\leq j\leq p$. 
Using this  identification we have that $\p^-\simeq M(q\times p,\C)$ and
$${\Cal D}=\left\{X\in M(q\times p, \C)\: \1_{p\times p} -X^*X>>0\right\}$$
and so by taking real points (i.e. $\sigma$-fixed points) we get
$$\Omega=\left\{ X\in M(q\times p, \R)\: \1_{p\times p} -X^tX>>0\right\}\ .\leqno(2.9)$$
\endexample

\head 3. The fine convex structure of $\Omega\simeq H/H\cap K$ \endhead 

The proof of the product formula for the $\c$-functions
depends on the fine convex geometry of $\Omega$ which we
discuss in this section. The problem we face here is that 
the $N^-$-part of $HP_{\text{\rm max}}\cap
\overline NP_{\text{\rm max}}=
\exp (\Omega )P_{\text{\rm max}}$ is not all of $N^-$. Hence the rank one reduction
depends on the understanding of the
the projection and intersection of $\Omega$ onto
subspaces corresponding to rank one subspaces corresponding
to each positive root. 
The results are taken from [16].

\ssk Define a convex cone in $\a$ by
$$c_{\text{\rm max}}=\left\{X\in\a\: (\forall \alpha\in\Delta_n^+) \ \alpha(X)\geq 0\right\}\ .$$
Then it follows readily from 2.8 that 
$$(\forall X\in c_{\text{\rm max}}) \qquad e^{\ad X}(\Omega)\subeq \Omega\ .\leqno(3.1)$$
This fact is closely related to the definition of the compression semigroup 
of $\Omega$ and the more geometric definition of causal symmetric spaces 
using cones (cf.\ Section 1). Later, in Section 5 when we will discuss spherical 
functions, we will entirely adopt this more geometric point view.   

\par For all $\alpha\in\Delta_n^-$ we let 
$$p_\alpha\: \n^-\to \g^\alpha$$
denote the projection along $\bigoplus_{\beta\in\Delta_n^-\atop\beta\neq \alpha} \g^\beta$. 
As mentioned earlier $\m=\z_{\k}(\a)=\z_\h(\a)$. Hence $N_H(\a)\subseteq N_K(\a )$ and so
$${\Cal W}_0\:=N_{H\cap K}/M=N_H(\a )/Z_H(\a)\ . $$
The group ${\Cal W}_0$ is called {\it the little Weyl group}. 
As $\Ad (k)T_0=T_0$ for all $k\in H\cap K$ it follows that $H\cap K$ normalizes
$\n^-$. In particular $w(\Delta_n^-)=\Delta_n^-$ for all $w\in {\Cal W}_0$.

\proclaim{Theorem 3.1} {\rm {\bf [16, Th. II.5]}} Let $\alpha\in \Delta_n^-$. Then
$$p_\alpha(\Omega)\subeq\Omega\ .$$
\endproclaim

\demo{Proof} A proof of Theorem 3.1 in full generality requires a fair amount of advanced structure theory 
which is beyond the scope of this exposition. Here we restrict ourselves for 
giving a proof for our basic example $\Sl(n,\R)/ \SO(p,q)$. The proof given below 
contains the idea of the general proof in [16]. 

\par Recall the facts and notations of Example 1.9 and Example 2.2. Let 
$X=\sum_{\alpha\in\Delta_n^-}X_\alpha\in\Omega$ with
$X_\alpha\in\g^\alpha$. We have to show that $X_\alpha\in\Omega$ for
all $\alpha\in\Delta_n^-$. Denoting by $\frak{S}_k$ the group of permutations
of $\left\{1,\ldots ,k\right\}$, then ${\Cal W}_0\simeq \frak{S}_q\times \frak{S}_p$ where the action
of a permutation $s=(s_1,s_2)$ is given by
$$s(\eps_{p+i}-\eps_j)=\eps_{p+s_1^{-1}(i)}-\eps_{s_2^{-1}(j)}\ , \quad
1\leq i\leq q, 1\leq j\leq p \ .$$
Hence all the roots in $\Delta_n^-$ are conjugate under the little Weyl group.
Thus we only have to show that
$X_\alpha\in\Omega$ for $\alpha=\eps_{p+1}-\eps_p$ (notice that $\Omega$ is invariant
under $N_{H\cap K}(\a)$).
If $q=1$ then there is only one root so there is nothing to prove. We can therefore assume
that $q>1$. Let $\delta={p-1\over q-1}$ and define  $Y\in\a$ by
$$Y=\diag(1, \ldots , 1, 0, 0, -\delta, \ldots , -\delta)$$
with $(p-1)$-times $+1$ on the diagonal. It is easy
to check that $Y\in c_{\text{\rm max}}$. Further notice that $\alpha(Y)=0$ while
$\beta(Y)<0$ for all $\beta\in \Delta_n^-$ with $\beta\neq \alpha$.
By  (3.1) we get
that
$$X_\alpha=\lim_{t\to \infty} e^{t\ad Y}( X) \in \Omega, $$
concluding the proof of the theorem. \quad\qed
\enddemo

\example{Example 3.2} {\bf (Main example continued)} For our main example $G/H=\Sl(n,\R)/ \SO(p,q)$
we use the natural identification $\n^-\simeq M(q\times p, \R)$.
Let
$$X=\sum_{1\leq i\leq q\atop 1\leq j\leq p} x_{ij} F_{ij}\in \n^-$$
with coordinates $x_{ij}\in \R$.
Recall that $X\in\Omega$ if and only if $\1_{p\times p} -X^tX>>0$. Thus if
$X\in\Omega$, then Theorem 3.1 implies that $x_{ij}F_{ij}\in\Omega$ for all $i,j$.
The latter is equivalent to 
$$\eqalign{&\1_{p\times p} -(x_{ij} F_{ij})^t (x_{ij} F_{ij})>>0\cr
\iff &\1_{p\times p} -x_{ij}^2  F_{jj}>>0\cr
\iff &|x_{ij}|< 1\cr}\ .$$
Thus if $X=\sum_{1\leq i\leq q\atop 1\leq j\leq p} x_{ij} F_{ij}\in \Omega$, then
$|x_{ij}|<1$ for all $i,j$. Notice that this is not obvious from 
(2.9).
\endexample

For $\alpha\in \Delta_n^-$  let $\g(\alpha)$ be the Lie algebra
generated by $\g^{\alpha}$ and $\g^{-\alpha}$. As $\tau (\g^{\alpha})=\theta (\g^{\alpha})
=\g^{-\alpha}$ it follows that
$\g(\alpha)$ is $\tau$ and $\theta$-stable. Hence $\tau_\alpha\:=\tau\res_{\g(\alpha)}$
and $\theta_\alpha\:=\theta\res_{\g(\alpha)}$ define involutions on $\g(\alpha)$. It is not hard to see
that $(\g(\alpha),\tau_\alpha)$ is a symmetric NCC pair  of rank one.
Write $\h(\alpha)=\h\cap \g(\alpha)=\g(\alpha)^{\tau_\alpha}$ and $\k(\alpha)=\k\cap\g(\alpha)=
\g(\alpha)^{\theta_\alpha}$. Let $G(\alpha)$,
$K(\alpha)$ be the
analytic subgroups of $G$ corresponding to $\g(\alpha)$ and $\k(\alpha)$. Set
$H(\alpha)=H\cap G(\alpha)=G(\alpha)^{\tau_\alpha}$. Finally
we let $\Omega(\alpha)\simeq H(\alpha)/H(\alpha)\cap K(\alpha)$ denote
the bounded symmetric domain associated to  $G(\alpha)/H(\alpha)$.
Then we obtain a canonic embedding of NCC-spaces
$$G(\alpha)/H(\alpha)\into G/H $$
and accordingly an embedding of the real bounded symmetric domain
$$\Omega(\alpha)\into \Omega\ .$$

\par The second main result on the fine convex structure of $\Omega$ is then:

\proclaim{Theorem 3.3} {\rm {\bf [16, Th. II.7]}} Let $\alpha\in \Delta_n^-$. Then
$$\Omega\cap \g^{\alpha}=\Omega(\alpha)\ . $$
\quad\qed
\endproclaim

\example{Example 3.4} {\bf (Main example continued)} Recall from
Example 2.2, that the
roots in $\Delta^-_n$ are $\eps_{p+i}-\eps_j$
with $1\leq i\leq q$ and $1\leq j\leq p$. In the
notation from that example we have
$$\Omega(\eps_{p+i}-\eps_j)=]-1,1[ F_{ij}\ .$$
Thus Theorem 3.3 means for our main example that 
$$\Omega\cap \R F_{ij}=]-1,1[ F_{ij}$$
for all $i,j$.
\endexample

\head 4. The $\c$-functions \endhead 

In this section we introduce the $\c$-function associated to the NCC space $G/H$. After
recalling its definition and basic properties  from [3] we  discuss the
product formula for the $\c$-function as a product of the $\c$-functions corresponding
to the rank one symmetric spaces $G(\alpha )/H(\alpha)$, $\alpha\in \Delta_n^+$.
This formula was established by J. Faraut in [2] for all CT spaces and
by P. Graczyk in [6] for our main example. The general formula
was finally obtained in [16]. The exposition here follows the geometric
approach {}from [16]; in particular the results of Section 3 will
enter in a crucial way.

\par Let $\Delta^+$ be any positive system in $\Delta$ containing $\Delta_n^+$.
Set $\Delta^-=-\Delta^+$ and define nilpotent subalgebras of $\g$ by
in the usual way by
$$\n=\bigoplus_{\alpha\in \Delta^+} \g^\alpha\quad\hbox{and}\quad
\oline \n=\bigoplus_{\alpha\in \Delta^-}\g^\alpha\ .$$
Then $\tau (\n)=\theta (\n)=\oline\n$.
Write $N$ and $\oline N$ for the analytic subgroups of $G$ corresponding to
$\n$ and $\oline \n$.

\par Recall that $HAN$ is open in $G$ and that the multiplication mapping
$$H\times A\times N\to G, \ \ (h,a,n)\mapsto han\leqno(4.1)$$
is a diffeomorphism onto $HAN$.  Write
$$HAN\ni x\mapsto (h_H(x),a_H(x),n_H(x))\in H\times A\times N\leqno(4.2)$$
for its invere
and notice the maps $h_H$, $a_H$, and  $n_H$ are analytic. For $x\in G$ we
also introduce the notation $x=k(x)a(x)n(x)$ with $k(x)\in K$, $a(x)\in A$,
and $n(x)\in N$.

\par For $\lambda\in\a_\C^*$ and  $a\in A$ we adopt the usual convention by setting
$a^\lambda=e^{\lambda (\log a)}$. Further we define
$\rho={1\over 2}\sum_{\alpha\in \Delta^+} (\dim \g^\alpha )\alpha$. 
\par Fix a Haar measure $d\oline n$ on $\oline N$  such that
$\int_{\oline N} a(\oline n)^{2\rho}\, d\oline n=1$.
We have now all ingredients to give a first definition of the subject proper
of this paper, the $\c$-function for $G/H$:

\definition {Definition 4.1} {\rm{\bf [3]}} Let $G/H$ be an NCC space. Then we
define {\it the $\c$-function of $G/H$ with parameter $\lambda\in \a_\C^*$}
by
$$\c(\lambda)=\int_{\oline N\cap HAN} a_H(\oline n)^{-(\lambda+\rho)} \ d\oline n\ ,\leqno(4.3)$$
whenever the integral converges absolutely. We denote by ${\Cal E}$ the
set of $\lambda\in \a_\C^*$ for which $\c(\lambda)$ is defined. 
\enddefinition 

The definition of the $\c$-function is motivated from the theory of spherical 
functions of $G/H$ developed in [3]  and which we will explain in the next section. 
Presently what is most disturbing  is the seemingly complicated structure 
of the domain of integration $\oline N\cap HAN$. Also it is not clear 
from the definition that ${\Cal E}\neq \eset$. 

\example{Example 4.2} Let us compute the $\c$-function for the rank one NCC
spaces starting with the simplest case of $G=\Sl (2,\R )$. Here 
$$H=\SO (1,1)=\pm \left\{h_t=\pmatrix \cosh t & \sinh t\\
\sinh t & \cosh t\endpmatrix \: t\in \R\right\}\ , $$
and
$$A= \left\{a_t=\pmatrix e^t & 0\\
0 & e^{-t}\endpmatrix \: t\in \R\right\}\ . $$
In this case $\Delta^+_n=\Delta^+$ and
$$N=N^+=\left\{n_x=\pmatrix 1 & x\\
0 & 1\endpmatrix \: x\in \R\right\}\qquad\hbox{and}\qquad
\overline N=N^-=\left\{\overline n_y=\pmatrix 1 & 0\\
y & 1\endpmatrix \: y\in \R\right\} . $$
A simple calculation now shows that
$$HP_{\text{\rm max}}=\left\{\pmatrix a & b\\
c & d\endpmatrix \: a^2-c^2>0\right\}\ ,  $$
$$\Omega\simeq ]-1,1[\ , $$
and
$$a_H\pmatrix a & b\\
c & d\endpmatrix =\pmatrix \sqrt{a^2-c^2}& 0\\
0 & 1/\sqrt{a^2-c^2}\endpmatrix \ . $$
Let us identify $\a_\C^*$ with $\C$ by $\lambda\mapsto \lambda \pmatrix 1 & 0\\
0 & -1\endpmatrix $. Then 
$$\eqalign{\c(\lambda)=\c_{\Omega}(\lambda)&=\frac{1}{\pi}\int_{-1}^1a(\overline n_y)^{\lambda +\rho}\, dx\cr
&= \frac{1}{\pi}\int_{-1}^1 (1-y^2)^{(\lambda+1)/2}\, dy\cr
&=\frac{1}{\pi}\int_{-1}^1t^{-1/2}(1-t)^{-(\lambda+1)/2}\, dt\cr
&=\frac{1}{\pi}B\left(\frac{1}{2},\frac{-\lambda+1}{2}\right)\ , \cr
}$$
where $B$ denotes the Beta function. 
\par For the general rank one case 
$G=\SO_o (n+1,1)/\SO_o(n,1)$ a very similar calculation applies, carried out 
in [3]. Here we also have $\Delta_n^+=\Delta^+=\left\{\alpha\right\}$. 
We then identify $\a_\C^*$ with
$\C$ by $\lambda \leftrightarrow (z/2)\alpha$. Then
$\rho =n$. Finally, there exists a constant $\kappa$ such
that
$$\c(\lambda)=\c_{\Omega}(\lambda )=\kappa B\left(\frac{n}{2},1-\frac{z+n}{2}\right)\ .$$
In particular,  it follows that $\c_\Omega$ extends to a meromorphic function
on all of $\a_\C^*$.
\endexample

\remark {Remark 4.3} The reader familiar with harmonic analysis on Riemannian 
symmetric spaces $G/K$ will recognize a similarity in the definitions
of the $\c$-function of $G/H$ with the one for $G/K$. In fact, replacing 
$H$ by $K$ in (4.3) gives the definition of the $\c$-function for $G/K$ (notice that 
$KAN=G$ so that $\oline N\cap KAN=\oline N$). This resemblance  is not by
coincidence. First of all both $\c$-functions show up as the coefficient of
the leading term in the assymptotic expansion of
spherical functions. The spherical function in question can be
written as an integral over $\overline N\cap HAN$ which results in the coefficient
beeing an integral over $\overline N\cap HAN$. We will discuss this in Section 5
(see also [22]). 
\par
A second reason, and more geometric, for this is
that the Riemannian symmetric spaces and the NCC spaces both belong 
to the  class of symmetric spaces $G/H$ which admit open $\Ad(H)$-invariant 
hyperbolic convex subsets in $\q$. This class of symmetric spaces 
(not necessary semisimple) was investigated geometrically in [12]
and subequently a unified approach to spherical functions was developed 
in [14,15]. 
\endremark

Let us now describe the geometry of the domain $\oline N\cap HAN$. 
For that we define $\Delta_0^\pm =\Delta^\pm\cap\Delta_0$ and 
set 
$$\n_0^+=\bigoplus_{\alpha\in\Delta_0^+}\g^\alpha\quad
\hbox{and}\quad \n_0^-=\bigoplus_{\alpha\in\Delta_0^-}\g^\alpha\ .$$
Then the triangular decomposition (2.1) implies 
$\Delta^+=\Delta_0^+\amalg \Delta_n^+$ as well as $\n=\n^+ \rtimes\n_0^+$
and $\oline \n=\n^-\rtimes \n_0^-$. Accordingly
if we let $N^\pm $ and $N_0^\pm$ be the analytic subgroups of 
$G$ corresponding to $\n^\pm,$ and $\n_0^\pm$, then we have 

$$N=N^+\rtimes N_0^+\quad\hbox{and}\quad \oline N=N^-\rtimes N_0^-\ .$$ 
\par The following proposition is well known (cf.\ [10, Ch.\ 5]). But because of 
its importance for us we will provide a quick proof.

\proclaim{Proposition 4.4} We have 
$$\oline N\cap HAN=\exp(\Omega)N_0^-\ .\leqno(4.4)$$
\endproclaim  
\demo{Proof} Let $P_{\text{\rm max}}=G(0)N^+$ as before. Recall that $M=Z_K(A)$ and define 
$$P_{\text{\rm min}}=MAN\ . $$ 
Then $P_{\text{\rm min}}$ is a minimal parabolic subgroup of $G$ contained in the maximal 
parabolic subgroup $P_{\text{\rm max}}$. By (2.7)
$HP_{\text{\rm max}}=\exp(\Omega) G(0) N^+$. As $HP_{\text{\rm max}}=HAN$ we get
$$\oline N\cap HAN =\oline N\cap \exp(\Omega)G(0)N_n^+\ . \leqno(4.5)$$
Let ${\Cal W}_0=N_{H\cap K}(\a)/M$ be as before the little
Weyl group and
${\Cal W}=N_{K}(\a)/M$ be the {\it big Weyl group}. Here we have
used that $M=Z_K(\a)=Z_H(\a)$.
Then we have the Bruhat decompositions

$$G(0)=\coprod_{w\in{\Cal W}_0} N_0^- wMAN_0^+w^{-1}\quad
\hbox{and}\quad G=\coprod_{w\in{\Cal W}} \oline N  wMANw^{-1}\ .\leqno(4.6)$$
Combining (4.5) and (4.6) now yields the assertion.
\quad\qed
\enddemo

\par As was observed in [3] the geometric splitting (4.4) of
the domain $\oline N \cap HAN$ implies  a splitting
of the $\c$-function. Define $\rho_0={1\over 2}
\sum_{\alpha\in\Delta_0^+} (\dim \g^\alpha) \alpha$. Normalize the
Haar measure $d\oline n_0$ on
$N_0$ such that $\int_{\oline N_0} a(\oline n_0)^{2\rho_0}\, d\oline n_0=1$.
Then the Haar measure $d\oline n=d\oline n_0d\oline n_n$ on
$\oline N$ satisfies $\int_Na(\oline n)^{2\rho}\, d\oline n=1$

\definition{Definition  4.5} For $\lambda\in \a_\C^*$ we define the partial
$\c$-functions
$$\c_0(\lambda)=\int_{N_0^-} a_H(\oline n_0)^{-(\lambda+\rho_0)} \ d\oline n_0$$
and
$$\c_\Omega(\lambda)=\int_{\exp(\Omega)} a_H(\oline n_n)^{-(\lambda+\rho)} \ d\oline n$$
whenever the defining integrals converge. We denote by ${\Cal E}_0$ and
${\Cal E}_\Omega$ the domains of definition of $\c_0$ and $\c_\Omega$, respectively.
\enddefinition

\par Notice that $\c_0(\lambda)$ is the $\c$-function of the Riemannian symmetric space
$G(0)/H\cap K$ and hence explicitely known ({\it Gindikin-Karpelevic formula} [5]).
{}From Definition 4.1,  (4.4) and the fact that $N^-$ is $G(0)$-invariant,
one  obtains the splitting
$$(\forall \lambda\in {\Cal E})\qquad \c(\lambda)=\c_0(\lambda)\c_\Omega(\lambda)\ .\leqno(4.7)$$
Thus in order to compute $\c(\lambda)$ it suffices to compute $\c_\Omega(\lambda)$. This was
accomplished by Faraut (cf.\ [2]) for CT spaces using Jordan algebra techniques in a very intelligent
way. Later, using a modification of Faraut's method, P. Graczyk
computed the $\c$-function for our main example
$\Sl (n,\R)/\SO (p,q)$ in [6]. But one of the
observations in [16] was that it is better to
compute $\c(\lambda)$ directly without using the
splitting (4.7). In fact the authors observed 
that all what was needed to make the Gindikin-Karpelevic idea work 
(cf.\ [5] and [4] for a particular nice exposition) in the NCC setup were the geometric 
results from Section 3. 

\par For  $\alpha\in \Delta$ let $H_\alpha\in [\g^\alpha ,\g^{-\alpha}]\subseteq\a$
be such that $\alpha (H_\alpha)=2$. Set
$m_\alpha=\dim \g^\alpha$.

\proclaim{Theorem 4.6} {\rm {\bf (The product formula for $\c_\Omega$) [16, Th. III.5]}} For the $\c$-function
$\c_\Omega$ of the real bounded symmetric domain $\Omega$ one has
$${\Cal E}_\Omega=\left\{ \lambda\in \a_\C^*\: (\forall \alpha\in \Delta_n^+)\
\Re \lambda(H_\alpha)< 2-m_\alpha\right\}$$
and
$$\c_\Omega(\lambda)=\kappa\prod_{\alpha\in\Delta_n^+} B\left({m_\alpha\over 2},
-{\lambda(H_\alpha)\over 2} -{m_\alpha\over 2} +1\right) $$
where $B$ denotes the Beta function and $\kappa$ is a positive constant only
depending on $(\g,\tau)$.\quad\qed
\endproclaim

\proclaim{Corollary 4.7} The function $\c_\Omega$ has a meromorphic continuation
to all of $\a_\C^*$.\quad\qed
\endproclaim

\example{Example 4.8} {\bf (Main example continued)} We write elements $\lambda\in \a_\C^*$ as
$\lambda=\sum_{j=1}^n \lambda_j \eps_j$ with coordinates $\lambda_j\in \C$.
Notice that $m_\alpha=1$ for all $\alpha\in \Delta$. Thus for $G/H=\Sl(n,\R)/\SO(p,q)$
we get
$${\Cal E}_\Omega=\left\{ \lambda\in \a_\C^*\: \Re (\lambda_j-\lambda_i)>-1, \ 1\leq i\leq p, \
p+1\leq j\leq n\right\}$$
and

$$\c_\Omega(\lambda)=\kappa \prod_{1\leq i\leq p\atop p+1\leq j\leq n} B\left({1\over 2}, {\lambda_j-\lambda_i\over 2}
+{1\over 2}\right)
=\kappa \pi^{pq\over2}\prod_{1\leq i\leq p\atop p+1\leq j\leq n} {\Gamma\left({\lambda_j-\lambda_i\over 2}
+{1\over 2}\right)\over \Gamma\left({\lambda_j-\lambda_i\over 2}
+1\right)}\ .\leqno(4.8)$$
\endexample

\remark {Problems 4.9} For the definition of the $\c$-functions one uses special positive systems
$\Delta^+$ which are compatible with the triangular decomposition (2.1). But what
happens if we allow arbitrary positive systems $\Delta^+$ of $\Delta$? Define
$N(\Delta^+)$ to be the analytic subgroup corresponding to
$$\n(\Delta^+)=\bigoplus_{\alpha\in \Delta^+}\g^\alpha$$
and $\overline N(\Delta^+):=\theta (N(\Delta^+))$.   Then we
can define a $\c$-function $\c_{\Delta^+}(\lambda )$ by
$$\c_{\Delta^+}(\lambda)=\int_{\overline N(\Delta^+)\cap HAN(\Delta^+)}
a_{H,\Delta^+}(\overline n)^{-(\lambda
+\rho(\Delta^+))}\, d\overline n$$
where $\rho (\Delta^+)$ and $a_{H,\Delta^+}$ are defined in the obvious way. 
Then some natural questions
are:
\item{$\bullet$} What is the geometry of $\oline N(\Delta^+)\cap HAN(\Delta^+)$?
\item{$\bullet$} For which $\lambda$'s - if any - does the integral defining
$\c_{\Delta^+}(\lambda )$ converge absolutely.
\item{$\bullet$} What is the significance of
$\c_{\Delta^+}(\lambda)$ in terms of harmonic analysis
on $G/H$ and in terms of representation theory?
\item{$\bullet$} Does a product formula for $\c_{\Delta^+}(\lambda)$ hold?
\endremark

\head 5. Spherical functions on $G/H$ \endhead 

In this section we recall a theory of spherical functions
for an NCC space $G/H$  developed 
in [3] and [19]. The connection with the
previous section is that the coefficient of the leading term in the asymptotic 
expansion of the spherical functions is given by the $\c$-function of $G/H$. 
The spherical functions introduced below will not be be defined on $G/H$ but on an open 
subset which is homemorphic to a hyperbolic cone $C\subeq \q$. 
In the sequel we will use the more geometric definition of NCC spaces through cones (cf.\ Definition 1.3)

 \par Recall the definiton of the maximal parabolic subgroup $P_{\text{\rm max}}=G(0)N$ of $G$. 
Write ${\Cal O}$ for the open $H$-orbit in the flag manifold $G/P_{\text{\rm max}}$ through the 
origin. {}From our discussion in Section 4 we have 

$${\Cal O}\simeq H/H\cap K\simeq \Omega\ .$$
We define the {\it compression semigroup} $S$ of ${\Cal O}$ by 

$$S=\left\{ g\in G\: g{\Cal O}\subeq {\Cal O}\right\}\ .$$
Clearly $S$ is an $H$-bi-invariant semigroup in $G$ contained in $HAN=HP_{\text{\rm max}}$. 
{}From (3.1) we obtain $\exp(c_{\text{\rm max}})\subeq S$.
If $S^0$ denotes the interior of $S$, then we have in addition
\ssk 
\item{$\bullet$} $S^0=H\exp(c_{\text{\rm max}}^0)H$ (cf.\ [10, Ch. 5]);
\item{$\bullet$} $C=\Ad(H)c_{\text{\rm max}}^0$ is an $\Ad(H)$-invariant open 
hyperbolic convex cone in $\q$ which contains no affine lines. Moreover 
the multiplication mapping 
$$H\times C\to S^0, \ \ (h,X)\mapsto h\exp(X)$$
is a homeomorphism (cf.\ [10, Th.\ 5.2.6]). 
\ssk

\example {Example 5.1} {\bf (Main example continued)} Even for $G/H=\Sl(n,\R)/ \SO(p,q)$
the cone $C$ forfeits a nice explicit description. However the trace of $C$ in $\p$ can 
be described nicely: 

$$\eqalign{C\cap\p=\Big\{ \pmatrix  A& 0\\  0 &B\endpmatrix 
\: & A\in M(p,\R),\  B\in M(q,\R);\ 
A^t=A,\  B^t=B,\  \tr (A) +\tr(B)=0, \cr 
& \min \Spec(A)>\max \Spec(B)\Big\}.\cr}$$
Notice that $T_0\in C$.
\endexample 

We fix a Haar measure $dh$ of $H$. 

\definition{Definition 5.2} Let $\lambda\in\a_\C^*$. Then we define the {\it spherical 
function with spectral parameter $\lambda$} by 
$$\phi_\lambda\: S^0\to\C, \ \ s\mapsto \int_H a_H(sh)^{\lambda-\rho}\ dh $$
whenever the defining integrals exist. 
\enddefinition

\remark {Remark 5.3} Several comments on the definition of $\phi_\lambda$ are appropriate: 
\par\nin (a) The functions $\phi_\lambda$ are only defined on $S^0$ and not on the 
whole symmetric space $G/H$. 

\par\nin (b) One can show that $\phi_\lambda$ is defined if and only if 
$\lambda\in {\Cal E}_\Omega$ (cf.\ [3, Th.\ 6.3]). 

\par\nin (c) All the spherical functions are $H$-bi-invariant and, as functions on $S^0$, 
common eigenfunctions of $\D(G/H)$, the algebra of $G$-invariant differential 
operators on $G/H$. 
\endremark

In order to relate the spherical functions with our $\c$-functions we need an integration 
formula whose proof can be found in [18].  Recall that $M=Z_H(\a)=Z_K(\a)$ is a compact subgroup of $H$.
Then we have a diffeomorphism 
$$N\cap HAN\to H/M, \ \ \oline n\mapsto h_H(\oline n) M\ .$$
Here $h_H(\oline n)\in H$ is defined by (4.2).

Write $m$ for a normalized Haar measure 
on $M$. Then we choose a Haar measure $d(hM)$ on $H/M$ such that 
$$\int_H f(h) \ dh =\int_{H/M} \int_{M} f(hm)\ dm \ d(hM)$$
holds for all $f\in L^1(H)$. 

\proclaim{Lemma 5.4} We can normalize $dh$  such that 
$$\int_{H/M} f(hM) \ d(hM) =\int_{\oline N \cap HAN} f(h_H(\oline n)M) a_H(\oline n)^{-2\rho}
\ d\oline n$$
holds for all $f\in L^1(H/M)$.\quad\qed 
\endproclaim

{}From Lemma 5.5 we now obtain an alternate definition 
for $\phi_\lambda$. 

\proclaim{Proposition 5.5} Let $\lambda\in {\Cal E}_\Omega$. 
Then for all $a\in S^0\cap A=\exp(c_{\text{\rm max}}^0)$ one has  
$$\phi_\lambda(a)=a^{\lambda-\rho}\int_{\oline N\cap HAN}
a_H(a\oline n a^{-1})^{-\rho+\lambda} a_H(\oline n)^{-(\lambda+\rho)}
\ d\oline n\ .$$ 
\endproclaim

\demo{Proof} Using Lemma 5.4 we obtain that 
$$\eqalign{\phi_\lambda(a)&=\int_H a_H(ah)^{\lambda-\rho} \ dh 
=\int_{H/M} a_H(ah)^{\lambda-\rho} \ d(hM)\cr 
&=\int_{\oline N\cap HAN} a_H(ah_H(\oline n))^{\lambda-\rho} 
a_H(\oline n)^{-2\rho}\  d\oline n\cr
&=\int_{\oline N\cap HAN} a_H(a\oline n a_H(\oline n)^{-1})^{\lambda-\rho} 
a_H(\oline n)^{-2\rho}\  d\oline n\cr
&=\int_{\oline N\cap HAN} a_H(a\oline n )^{\lambda-\rho} 
a_H(\oline n)^{-(\lambda+\rho)}\  d\oline n\cr
&=a^{\rho-\lambda}\int_{\oline N\cap HAN} a_H(a\oline n a^{-1} )^{\lambda-\rho} 
a_H(\oline n)^{-(\lambda+\rho)}\  d\oline n\ , \cr}$$
as was to be shown.
\quad\qed 
\enddemo

Denote by $\phi_\lambda^0$ the spherical function
$$\phi_{\lambda}^0(g)=\int_{K\cap H}a(gk_0)^{\lambda -\rho_0}\ dk_0\qquad (g\in G(0)) $$
for the Riemannian symmetric space $G(0)/K\cap H$.  Here
$\rho_0=\frac{1}{2}\sum_{\Delta_0^+}m_\alpha \a$ as before and $dk_0$ denotes
the normalized Haar measure on $K\cap H$. Finally let $\rho_+=\rho-\rho_0
=\sum_{\alpha\in \Delta^+_n}m_\alpha \alpha$.
By Proposition 5.5 we get:

\proclaim {Lemma 5.6} Let $s\in S^0$ and $\lambda \in {\Cal E}_0\cap {\Cal E}_{\Omega}$.
Then
$$\lim_{t\to \infty}e^{t(\rho-\lambda )(T_0)}\phi_\lambda (\exp (tT_0)s)
=\c_{\Omega}(\lambda)\phi^0_{\lambda-\rho_+} (s)\ . $$
In particular, $\phi_{s\lambda}=\phi_\lambda$ for all $s\in {\Cal W}_0$.\quad\qed 
\endproclaim

Define an open Weyl chamber $A^+\subseteq A$ by 
$$A^+=\left\{a\in A\: (\forall \alpha\in
\Delta^+)\ a^\alpha >1\right\}=\left\{ a\in A\: (\forall \alpha\in \Delta^+)\ \alpha(\log a)>0\right\}\ .$$
Notice that $A^+\subeq S^0\cap A$. In fact
$S^0\cap A=(\bigcup_{w\in {\Cal W}_0}\overline{wA^+w^{-1}})^0$. 
Then we obtain from Proposition 5.5: 

\proclaim{Proposition 5.7}  Let $\lambda\in {\Cal E}_\Omega$. Then
$$\lim_{a\to\infty\atop a\in A^+} a^{\rho-\lambda}\phi_\lambda(a)=\c(\lambda)\ .$$
\endproclaim

\demo{Proof} Notice that $\lim_{A^+\ni a\to \infty} a(a\oline na^{-1})=\1$
for all $\oline n\in \oline N$. Thus the assertion  follows from Proposition 
5.5 provided we can interchange the limit with the integral sign. 
In fact, one can show that dominated convergence applies and we refer to 
the proof of [3, Th. 6.8] for the technical details. 
\quad\qed 
\enddemo

It was shown in [19]  how to obtain from  Lemma 5.6 and Proposition 5.7 
a Harish-Chandra type expansion formula for the spherical function $\phi_\lambda$.
To be more precise let $\Lambda = \N_{0}[ \Delta^+]$. Following Harish-Chandra
we define for $\mu\in \Lambda$
a rational function $\a_\C^*\ni \lambda\mapsto\Gamma_\mu (\lambda)\in \C$ by
$\Gamma_0(\lambda)=1$ and then
$$((\mu ,\mu )-2(\mu ,\lambda ))\Gamma_\mu (\lambda)=2\sum_{\alpha\in\Delta^+}
m_\alpha \sum_{k\in \N}\Gamma_{\mu -2k\alpha}
(\lambda)(\mu+\rho-2k\alpha+\lambda , \alpha )\ . $$
Then $\Gamma_\mu$ is holomorphic as a function of $\lambda$ in
a open dense subset of $\a_\C^*$. We say that $\lambda$ is {\it generic}
if $\Gamma_\mu$ is holomorphic in the point $\lambda$. 
Furthermore we define functions $\Phi_\lambda$ on  $A^+$ by
$$\Phi_\lambda (a)=a^{\lambda-\rho}\sum_{\mu \in \Lambda}\Gamma_\mu(\lambda)
a^{-\mu}\ .$$ 
With these notations and supressing further details we then have:

\proclaim{Proposition 5.8} {\bf [19](Expansion formula)} For $\lambda\in
{\Cal E}_{\Omega}$ generic and $a\in A^+$ we have
$$\varphi_\lambda(a)=\sum_{w\in {\Cal W}_0}\c(w\lambda)\Phi_{w\lambda}(a)\ .$$
\quad\qed
\endproclaim

\proclaim{Corollary 5.9} Let $a\in A^+$. Then the function $\lambda\mapsto \varphi_\lambda (a)$
extends to a meromorphic function on $\a_\C^*$.\quad\qed
\endproclaim

\par Knowing that $\phi_\lambda$ extends to a meromorphic function on $\a_\C^*$ 
in the parameter $\lambda$ it is natural to study the singularities of this function. 
It is also not clear from the above that these meromorphically 
extended functions can still be
defined as $H$-invariant functions on $S^0$. Both of these problems were
addressed in [23] and [26].

\proclaim{Theorem 5.10} {\bf ([23, Cor.\ 8.2])} Denote by $n_\Omega (\lambda)$
the numerator of $\c_\Omega (\lambda)$.
Then there exists a ${\Cal W}_0$-invariant tubular neighborhood $U\subseteq A_\C$ of
$S^0\cap A$ such that
$$(\lambda ,a)\mapsto \frac{\phi_\lambda(a)}{n_\Omega (\lambda )}$$
extends to a ${\Cal W}_0\times{\Cal W}_0$-invariant holomorphic function
on $\a_\C^*\times U$.\quad\qed
\endproclaim
 
\head 6. Spherical highest weight representations \endhead 

In the last two remaining sections we discuss the relevance of the $\c$-function for 
representation theory. In this section we will explain the relation 
of $\c(\lambda)$ with the formal dimensions of representations of the
holomorphic discrete series on the 
$c$-dual space $G^c/H$. For a fixed choice of $\Delta^+$ the
representations of the holomorphic discrete series are by definition 
$H$-spherical highest weight representations
that can be realized as discrete summand in $L^2(G/H)$ (cf. [24]).
We will begin our discussion with some brief recall on spherical representations
and the $H$-spherical distribution vectors. Then we explain the
connection between the $\c$-function and the formal dimension  of
spherical representations.

\subhead General facts about spherical representations \endsubhead 

\par In this subsection $G$ will denote an arbitrary Lie group, $H$ 
a symmetric subgroup of $G$, and 
$(\pi, {\Cal H})$ a unitary representation of $G$. We write 
${\Cal H}^\infty$ and ${\Cal H}^\omega$ for the $G$-modules of smooth respectively analytic 
vectors. The strong antiduals of ${\Cal H}^\infty$ and ${\Cal H}^\omega$ are denoted by 
${\Cal H}^{-\infty}$ and ${\Cal H}^{-\omega}$ and refered to as the $G$-modules of 
{\it distribution} respectively {\it hyperfunction  vectors} of $(\pi, {\Cal H})$. 
Notice the following sequence of $G$-invariant continuous inclusions
$${\Cal H}^\omega \into {\Cal H}^\infty \into {\Cal H} \into {\Cal H}^{-\infty} \into {\Cal H}^{-\omega}\ .$$ 
If $G$ is reductive and $K<G$ is a maximal compact subgroup,
then ${\Cal H}_K$ will denote the $(\g,K)$-module
of $K$-finite vectors of $(\pi, {\Cal H})$.
Notice that if $(\pi ,{\Cal H})$ is irreducible then ${\Cal H}_K\subseteq {\Cal H}^{-\omega}$.

\par If $L$ is a group and $V$ is an $L$-module, then we denote by $V^L$ the space of $L$-fixed vectors 
on $V$. Similarily if $\l$ is a Lie algebra
and $V$ is a $\l$ module,
then we set $V^\l=\left\{v\in V\: (\forall X\in \l)\, Xv=0\right\}$.  
Finally if $V$ is a vector space, then we denote by $V^{\oline *}$ its algebraic antidual. 

\par Let $L<G$ be a subgroup. Then an irreducible unitary representation 
$(\pi, {\Cal H})$ of $G$ is called {\it $L$-spherical} if $({\Cal H}^{-\omega})^L\neq \left\{0\right\}$.
A nonzero element in $({\Cal H}^{-\omega})^L$ is called a {\it $L$-spherical
hyperfunction vextor } and similarly we call a nonzero element in
$({\Cal H}^{-\infty})^L$ a {\it $L$-spherical distribution vector}.
In general it is easier to work with
the $L$-spherical hyperfunction vectors, but the
following theorem of Brylinski, van den Ban, and Delorme (for a proof see [1]),
clarifies the relation between those two notions in
our case. 

\proclaim{Theorem 6.1} Suppose that $G$ is reductive and $H<G$ is a symmetric subgroup, then 
we have $({\Cal H}^{-\omega})^H=({\Cal H}^{-\infty})^H$. 
If in addition $H$ is connected, then
$$({\Cal H}^{-\omega})^H=({\Cal H}^{-\infty})^H =({\Cal H}_K)^{\oline *, \h}\, .$$
\quad \qed 
\endproclaim

\subhead Spherical highest weight modules \endsubhead

{}From now on $(\pi, {\Cal H})$ will denote a unitary highest weight representation of 
the Hermitian group  $G^c$. Recall that this means that there exists an irreducible 
$K^c$-submodule $V\subeq {\Cal H}$ such that 

\item{} {\bf (HW1)} ${\Cal H}_{K^c}={\Cal U}(\g_\C)V$;
\item{}{\bf (HW2)} $\p^+ V=\left\{0\right\}$.

\par\nin Here ${\Cal U}(\g_\C )$ denotes the {\it universal enveloping algebra}
of $\g_\C$. From now on we will also assume that $(\pi, {\Cal H})$ is $H$-spherical. This
implies in particular that the  $K^c$-module $V$ above  is $H\cap K^c$-spherical. Denote by
$v_0$ a nontrivial element in $V^{H\cap K^c}$.

\par The next step is to show that $H$-spherical highest weight modules have multiplicity
one, i.e. $\dim_\C ({\Cal H}^{-\omega})^H=1$. This
was first done in [24] for the holomorhic discrete series, and
the general case was done in [9].  The next lemma contains 
the crucial information.  

\proclaim{Lemma 6.2} Let $(\pi, {\Cal H})$ be an $H$-spherical unitary highest
weight representation of 
$G^c$. Then, as a $\h$-module,  ${\Cal H}_{K^c}$ is a quotient of ${\Cal U}(\h)\otimes_{{\Cal U}(\h\cap \k^c)}
V$. 
\endproclaim

\demo{Proof} Following the proof of Lemma 3.1.1 in [9] we note
that $\g_\C^c=\g_\C=\h_\C +\a_\C +\n_\C$ and so 
${\Cal U}(\g_\C^c)={\Cal U}(\h_\C) {\Cal U}(\a_\C +\n_\C)$. 
It follows from (HW2) and the fact that $V$ is a $\k^c$-module
that $(\a_\C +\n_\C)V\subeq V$. 
Thus ${\Cal H}_{K^c}={\Cal U}(\h_\C) V$ from which the lemma easily follows. 
\quad\qed 
\enddemo

\proclaim{Proposition 6.3} {\rm \bf (Multiplicity one)} Let $(\pi, {\Cal H})$
be an $H$-spherical unitary 
highest weight representation of $G^c$. Then $\dim ({\Cal H}^{-\omega})^H=1$. 
\endproclaim 

\demo{Proof} We follow the proof of Lemma 3.1.2 in [9]. By Theorem 6.1 it is enough to
show that $\dim ({\Cal H}_{K^c}^{\oline *})^\h= 1$. According 
to Lemma 6.2 this will be implied by  $\dim [{\Cal U}(\h)
\otimes_{{\Cal U}(\h\cap \k^c)}V ]^{\oline *, \h} =1$. 
for that write ${\Cal U}(\h)_+$ for the subspace of ${\Cal U}(\h)$  
consisting of elements with zero constant term. Then any $\h$-invariant 
functional on ${\Cal U}(\h)\otimes_{{\Cal U}(\h\cap \k^c)}V$ must necessarily vanish on 
${\Cal U}(\h)_+\otimes_{{\Cal U}(\h\cap \k^c)}V$. Thus up to scalar $\nu$ is given by 

$$\eqalign{\nu(h\otimes v)&= 0 \qquad \text{\rm if}\  h\in {\Cal U}(\h)_+\cr 
 \nu(h\otimes v)&= \la v_0, v\ra \qquad \text{\rm if}\  h=\1\ .}\leqno (6.1) $$
This proves the proposition. \quad\qed 
\enddemo
 
Write $\lambda$ for the highest weight of the $K^c$-module $V$. As $V$ is $H\cap K^c$-spherical
we have $\lambda\in \a^*$. In the sequel we write $(\pi_\lambda, {\Cal H}_\lambda)$
and $V_\lambda$ instead of $(\pi, {\Cal H})$ and $V$ in order to indicate the dependance on 
$\lambda$. According to Proposition 6.3 we have 
$$({\Cal H}_\lambda^{-\omega})=\C \nu$$
is one dimensional. For the appropriate normalization
of $\nu$ let $v_\lambda$ be a normalized highest weight vector 
of ${\Cal H}_\lambda$. As $\la v_0, v_\lambda\ra \neq 0$ it follows from (6.1)
that $\la \nu, v_\lambda\ra \neq 0$. We can therefore normalize 
$\nu$ by
$$ \la \nu, v_\lambda\ra =1.\leqno{(6.2)\ \hbox{(Normalization)}}$$

\remark{Remark 6.4} (a) The motivation for the normalization (6.2) comes from the group case 
$G=H\times H$. In this case all $H$-spherical
highest weight representations of $G$ are of the form 
$(\pi\otimes\overline \pi, {\Cal H}\hat \otimes \overline{\Cal H})$
where $(\pi, {\Cal H})$ a highest weight representation 
of $H$ and    $(\overline\pi, \overline{\Cal H})$ is its dual representation. 
\par Write $B_1({\Cal H})$ and $B_2({\Cal H})$ for the
spaces of trace class and Hilbert-Schmidt operators 
on ${\Cal H}$. We will identify ${\Cal H}\hat \otimes \overline{\Cal H}$ with  $B_2({\Cal H})$. 
Then $B_2({\Cal H})^{\infty}\subeq B_1({\Cal H})$ (cf.\ [8, Th.\ A.2]) and the conjugate trace 
$$\overline{\tr} \: B_2({\Cal H})^\infty\to \C, \ \ X\mapsto \oline {\tr X}$$
is a canonically normalized $H$-invariant distribution vector. It is easy to see that 
$\overline{\tr}$ satisfies (6.2).
\par\nin (b) Unitary spherical highest weight representations are meanwhile quite well 
understood. For more and an almost complete classification let us refer to [13].
\endremark

The next theorem is the main result in this section. It tells us how to obtain $\nu$ from 
$v_\lambda$. 

\proclaim{Theorem 6.5} {\rm \bf [11, Th. 2.16] (Averaging Theorem)} Let $(\pi_\lambda, {\Cal H}_\lambda)$ be an 
$H$-spherical unitary highest weight representation with paramter $\lambda$ 
satisfying $\lambda+\rho \in {\Cal E}_\Omega$. Then we have 
$$\int_H \pi_\lambda(h)v_\lambda \ dh =\c(\lambda+\rho)\nu\ .\leqno(6.3)$$
Here the left hand side is a convergent ${\Cal H}_\lambda^{-\omega}$-valued integral. 
\endproclaim

\demo{Proof} (Sketch) We will not go into the technical details of the convergence 
of the left hand side (LHS) in (6.3) and refer to the proof of [11, Th.\ 2.16] instead. 
So let us assume that the LHS actually converges. Because of  ``multiplicity 
one'' (cf.\ Proposition 6.3),  the LHS must be a multiple of $\nu$. In view of our 
normalization of $\nu$, it 
is hence enough to check that 
$$\int_H \la \pi_\lambda(h)v_\lambda, v_\lambda\ra \ dh =\c(\lambda+\rho)\ .$$
Using the integration formula from Lemma 5.4 we obtain that 
(note that $v_\lambda$ is $M$-fixed)
$$\eqalign{\int_H \la \pi_\lambda(h)v_\lambda, v_\lambda\ra \ dh &=
\int_{H/M} \la \pi_\lambda(h)v_\lambda, v_\lambda\ra \ d(hM)\cr 
&= \int_{\oline N \cap HAN}  \la \pi_\lambda(h_H(\oline n))v_\lambda, v_\lambda\ra  
a(\oline  n)^{-2\rho} \ d\oline n \ . \cr}$$
As $\oline n=h(\oline n)a_H(\oline n)n_H(\oline n)$ we
get $h_H(\oline n)=\oline n  a_H(\oline n)^{-1} n $ for some $n\in N$. 
Write $\pi_\lambda^V$ for the representation of $K_\C^c$ on $V_\lambda$. 
Using this and  the fact that $d\pi_\lambda(X)v_\lambda=0$ for all $X\in \n$ we
obtain for  the integrand that 
$$\la \pi_\lambda(h_H(\oline n))v_\lambda, v_\lambda\ra=
\la \pi_\lambda^V(a_H(\oline n^{-1}))v_\lambda, v_\lambda\ra =a_H(\oline n)^{-\lambda}\ .$$
Finally we get 
$$\int_H \la \pi_\lambda(h)v_\lambda, v_\lambda\ra \ dh 
=\int_{\oline N \cap HAN}  a_H(\oline n)^{-(\lambda+2\rho)} d\oline n=\c(\lambda+\rho)\ , $$
as was to be shown. \quad\qed 
\enddemo

\remark{Remark 6.6} {\bf (Formal dimension)} Theorem 6.5 features many applications. 
One of them is the determination of the formal dimension of a representation 
$(\pi_\lambda, {\Cal H}_\lambda)$ of the holomorphic discrete series 
on $G^c/H$. 
For such a representation we have a $G^c$-equivariant embedding 
$${\Cal H}_\lambda^\infty\to C^\infty (G^c/H), \ \ v\mapsto 
f_v; \ f_v(gH)=\la \pi_\lambda(g^{-1})v, \nu\ra$$
which extends to a continuous $G$-equivariant map 
${\Cal H}_\lambda\to L^2(G^c/ H)$ satisfying
$$(\forall v,w\in {\Cal H}_\lambda)\qquad {1\over d(\lambda)} \la v,w\ra 
=\la f_v, f_w\ra. \leqno(6.4)$$ 
Here $d(\lambda)>0$ is a positive number refered to as the {\it formal dimension}. 
Then one can essentially deduce from (6.3) the formula
$$d(\lambda)=\c_{G/H}(\lambda+\rho_H) \c_{G_\C/G}(\lambda+\rho_G)^{-1}\leqno (6.5)$$  
(cf.\ [11, Th.\ 3.6, Th.\  4.15, Cor. 4.16]). In (6.5) the functions $\c_{G/H}$ and 
$\c_{G_\C/ G}$ stand for the $\c$-functions of the NCC spaces $G/H$ and $G_\C/G$ (cf.\ Example 1.8); 
further $\rho_H$ and $\rho_G$ refer ro the different $\rho$'s of  the NCC pairs 
$(\g,\h)$ and $(\g_\C, \h_\C)$.  
\par In the group case $G=H\times H$ the formula (6.4) simplifies to 
$$d(\lambda)=\c_{H_\C/H}(\lambda+\rho)^{-1} \leqno(6.5)$$
and Harish-Chandra's formula for the formal dimension for the holomorphic discrete 
series on $G/H\simeq H$ (cf.\ [7])
$$d(\lambda)=\prod_{\alpha\in \Delta^+} {\la \lambda+\rho, \alpha\ra \over \la \rho,\alpha\ra}
\leqno(6.6)$$
becomes a special case of the product formula in Theorem 4.6 and the Gindikin-Karpelevic 
formula. 
\endremark

\head 7. Spherical characters of spherical highest weight representations \endhead 

In this section we will give another application of the Averaging Theorem 6.5.
We  relate the $H$-spherical distribution vector $\nu$ to
the spherical function  $\varphi_{\lambda +\rho}$ and the
normalizing factor turns out to be exactly the $\c$-function. Combined
with the expansion formula in Proposition 5.8 this gives  a
generalization of the Harish-Chandra character formula for the 
holomorphic discrete series  on a compact
Weyl chamber. For a different approach for the character formula using the
Hardy space realization of the spherical highest weight represenation
in the holomorphic discrete series see [20,21]. 

\par Throughout this section $(\pi_\lambda, {\Cal H}_\lambda)$ denotes   
an $H$-spherical unitary highest weight representation of $G^c$. The
corresponding
{\it spherical character} is given by 

$$\theta_\lambda\: C_c^\infty(G^c/H)\to \C, \ \ f\mapsto\la \pi_\lambda(f)\nu,\nu\ra
=\la \int_{G^c/H} f(gH)\pi(g)\nu \ dgH, \nu\ra\ .$$  
It is a general fact that $\theta_\lambda$ defines an $H$-invariant
distribution on $G^c/H$. 
\par For the group case a deep result of Harish-Chandra asserts that $\theta_\lambda$
is given by a locally integrable function. For general symmetric spaces this $G/H$ 
this becomes unfortuntately false. For a nice
formula for $\theta_\lambda$ one needs to extend $\theta_\lambda$
to a holomorphic function in a complex domain in $G_\C/H_\C$, i.e., concider
$\theta_\lambda$ as a hyperfunction. 
\par In this section we will therefore give a short discussion
of the analytic continuation of $\theta_\lambda$
along  $H\exp(i\a)H/H$ in $G^c/H$. Firstly,  it follows
from a beautiful result of Olshanski and Stanton that
the representation $\pi_\lambda$ can be 
analytically continued to a representation of the open semigroup $S^0=H\exp (c_{\text{\rm max}}^0)H$
(cf.\ [25], [27]).  
Moreover, by a result of K.-H. Neeb [17], the operators
$\pi_\lambda(s)$ for $s\in S^0$ are all of trace 
class and by [14, Appendix] we have 

$$(\forall s\in S^0) \qquad \pi_\lambda(s){\Cal H}_\lambda^{-\omega}
\subeq {\Cal H}_\lambda^\omega\ . \leqno(7.1)$$ 
 
{}From this it follows that
$$\Theta_\lambda\: S^0\to\C, \ \ s\mapsto \la \pi_\lambda(s)\nu, \nu\ra$$
defines an $H$-invariant analytic function on $S^0$.
We call the parameter $\lambda$ {\it regular} if ${\Cal H}_{\lambda, K^c}$
is isomorphic to ${\Cal U}(\g_\C)\otimes_{{\Cal U}(\k_\C +\p^+)} V_\lambda$
as $\g_\C$-module. Notice that $\lambda$ is regular for all parameters
belonging to the holomorphic discrete series (cf.\ [13, Th.\ A.1]).

\proclaim{Theorem 7.1} {\rm \bf(Spherical Character formula)}. Let $(\pi_\lambda, {\Cal H}_\lambda)$
be an $H$-spherical unitary highest weight representation of $G^c$ with regular 
parameter $\lambda$. Then we have 
$$(\forall s\in S^0)\qquad \Theta_\lambda(s)={1\over \c(\lambda+\rho)}\cdot \phi_{\lambda+\rho}(s)\ .$$
\endproclaim

\demo{Proof} We reproduce the proof from [11, Th.\ 5.4] for all $\lambda$ with $\lambda+\rho\in {\Cal E}_\Omega$. 
The analytic continuation to all regular parameters was established in [9, Th.\ 4.1.2] using 
results from [19]. 
So let us assume that $\lambda+\rho\in {\Cal E}_\Omega$. Then it follows from Theorem 6.5 that 

$$\eqalign{\Theta_\lambda(s)&=\la \pi_\lambda(s)\nu, \nu\ra\cr 
&={1\over \c(\lambda+\rho)} \int_H \la \pi_\lambda(sh)v_\lambda, \nu\ra \ dh \cr 
&={1\over \c(\lambda+\rho)} \int_H a_H(sh)^{\lambda}\ dh\cr 
&= {1\over \c(\lambda+\rho)}\cdot \phi_{\lambda+\rho}(s), \cr }
$$
as was to be shown. \quad\qed  
\enddemo

\refstyle{A}

\enddocument